\documentclass[12pt]{amsart}
\usepackage{amscd,amssymb,amsmath,multicol,multirow}
\usepackage{float}
\restylefloat{table}
\newtheorem{thm}[equation]{Theorem}
\numberwithin{equation}{section}

\begin{document}
\raggedbottom \voffset=-.7truein \hoffset=0truein \vsize=8truein
\hsize=6truein \textheight=8truein \textwidth=6truein
\baselineskip=18truept
\def\mc{\multicolumn}
\def\ss{\smallskip}
\def\ssum{\sum\limits}
\def\dsum{{\displaystyle{\sum}}}
\def\la{\langle}
\def\ra{\rangle}
\def\on{\operatorname}
\def\gcd{\on{gcd}}
\def\a{\alpha}
\def\bz{{\Bbb Z}}
\def\eps{\epsilon}
\def\br{{\bold R}}
\def\bc{{\bold C}}
\def\bN{{\bold N}}
\def\nut{\widetilde{\nu}}
\def\tfrac{\textstyle\frac}
\def\b{\beta}
\def\G{\Gamma}
\def\g{\gamma}
\def\zt{{\bold Z}_2}
\def\zth{{\bold Z}_2^\wedge}
\def\bs{{\bold s}}
\def\bg{{\bold g}}
\def\bof{{\bold f}}
\def\bq{{\bold Q}}
\def\be{{\bold e}}
\def\line{\rule{.6in}{.6pt}}
\def\xb{{\overline x}}
\def\xbar{{\overline x}}
\def\ybar{{\overline y}}
\def\zbar{{\overline z}}
\def\ebar{{\overline \be}}
\def\nbar{{\overline n}}
\def\rbar{{\overline r}}
\def\Ubar{{\overline U}}
\def\et{{\widetilde e}}
\def\ni{\noindent}
\def\ms{\medskip}
\def\ehat{{\hat e}}
\def\xhat{{\widehat x}}
\def\nbar{{\overline{n}}}
\def\minp{\min\nolimits'}
\def\N{{\Bbb N}}
\def\Z{{\Bbb Z}}
\def\S{{\Bbb S}}
\def\M{{\Bbb M}}
\def\C{{\Bbb C}}
\def\el{\ell}
\def\mo{\on{mod}}
\def\lcm{\on{lcm}}
\def\dstyle{\displaystyle}
\def\Remark{\noindent{\it  Remark}}
\title
{Enumerating lattices of subsets}
\author{Donald M. Davis}
\address{Department of Mathematics, Lehigh University\\Bethlehem, PA 18015, USA}
\email{dmd1@lehigh.edu}
\date{December 5, 2013}

\maketitle
\begin{abstract} If $X_1,\ldots,X_k$ are sets such that no one is contained in another, there is an associated lattice on $2^{[k]}$
corresponding to inclusion relations among unions of the sets. Two lattices on $2^{[k]}$ are equivalent
if there is a permutation of $[k]$ under which they correspond. We show that for $k=1$, 2, 3, and 4, there are 1, 1, 4, and 50 equivalence classes
of lattices on $2^{[k]}$ obtained from sets in this way.
We cannot find a reference to previous work on this enumeration problem in the literature, and so wish to introduce it for subsequent investigation.
We explain how the problem arose from algebraic topology.
\end{abstract}
\section{Introduction}
For a possible application to algebraic topology, we have become interested in an enumeration problem
for lattices of subsets, which we have been unable to find in the literature. We wish  to introduce it for further investigation.

Let $[k]=\{1,\ldots,k\}$, and $2^{[k]}$ its power set.
If $\M=\{X_1,\ldots,X_k\}$ is a collection of sets, and $S\subset [k]$, let \begin{equation}\label{Meq}\M_S:=\bigcup_{i\in S}X_i.\end{equation}
We say that $\M$ is {\it proper} if it is never the case that $X_i\subset X_j$ for $i\ne j$.
Any $\M$ defines a lattice $L(\M)$ on $2^{[k]}$ by $S\le T$ if $\M_S\subset \M_{T}$.
Lattices $L$ and $L'$ on $2^{[k]}$ are said to be equivalent if there is a permutation $\sigma$ of $[k]$ under which the induced permutation of $2^{[k]}$
preserves the lattice relations; i.e., $\sigma(S)\le' \sigma(T)$ iff $S\le T$.
We wish to enumerate the equivalence classes of all possible $L(\M)$'s for proper $\M$'s of size $k$.

For $k=1$ and $k=2$, there is only one equivalence class each of such $L(\M)$'s. Indeed, for $k=2$, we must have $\M_{\{1\}}$ and $\M_{\{2\}}$ both strictly contained in $\M_{\{1,2\}}$.
For $k=3$, there are four equivalence classes. Out of the three possible nontrivial inclusions $\M_{\{1\}}\subset \M_{\{2,3\}}$, $\M_{\{2\}}\subset \M_{\{1,3\}}$, and $\M_{\{3\}}\subset\M_{\{1,2\}}$,
it can be the case that 0, 1, 2 or 3 of these hold, giving the four equivalence classes. We  also observe
that each of the four equivalence classes can be realized by actual sets. See Table \ref{t1} for realizations.

\begin{table}[H]
\caption{$k=3$}\label{t1}
\begin{tabular}{c|ccc}
Inclusions&$X_1$&$X_2$&$X_3$\\
\hline
$\emptyset$&$\{a\}$&$\{b\}$&$\{c\}$\\
$\M_{\{3\}}\subset\M_{\{1,2\}}$&$\{a,c\}$&$\{b,d\}$&$\{c,d\}$\\
$\M_{\{3\}}\subset\M_{\{1,2\}},\ M_{\{2\}}\subset\M_{\{1,3\}}$&$\{a,b,d\}$&$\{b,c\}$&$\{c,d\}$\\
$\M_{\{3\}}\subset\M_{\{1,2\}},\ M_{\{2\}}\subset\M_{\{1,3\}},\ \M_{\{1\}}\subset\M_{\{2,3\}}$&$\{a,b\}$&$\{b,c\}$&$\{c,a\}$
\end{tabular}
\end{table}

Note that the condition that $\M$ be proper forces that these inclusions are strict, and each of these inclusions
is equivalent to a corresponding equality involving $\M_{\{1,2,3\}}$. For example, $\M_{\{1\}}\subset \M_{\{2,3\}}$ iff $\M_{\{2,3\}}= \M_{\{1,2,3\}}$.
For consideration of equivalence classes of $L(\M)$'s when $k>3$, we find the equality viewpoint to be more convenient than the containment viewpoint, which leads
to the following alternative formulation.

Let $\S_k$ denote the collection of subsets of $[k]$ of cardinality greater than 1. A {\it configuration of size} $k$ is an equivalence relation on $\S_k$ such that
$S\sim T$ implies $S\cup\{i\}\sim T\cup\{i\}$ for all $i$. Configurations $\sim$ and $\sim'$ are equivalent if there is a permutation $\sigma$ of $[k]$ such that $\sigma(S)\sim'\sigma(T)$ iff $S\sim T$. A configuration $\sim$ is {\it realizable} if there exists a proper  collection $\M$  as above
such that $\M_S=\M_T$ iff $S\sim T$.

Note that we do not include singleton sets in $\S_k$ because of the ``proper'' condition on $\M$. We would never have $\M_{\{i\}}=\M_T$ if $|T|>1$. So if we included
singleton sets, they would be only equivalent to themselves.

The four equivalence classes of configurations of size 3 are those in which 0, 1, 2, or 3 of the sets $\{1,2\}$, $\{1,3\}$, and $\{2,3\}$ are $\sim\{1,2,3\}$.

Our main theorem is
\begin{thm} There are exactly $50$  equivalence classes of configurations of size $4$, and each is realizable.\end{thm}

We felt that this enumeration problem would be of the type that would have been studied and have its results appearing on Sloane's
website (\cite{Sloane}), but the sequence 1, 1, 4, 50 does not appear there in any context related to counting sets.
Nor have we been able to find references to this specific problem in the literature.

This enumeration problem, in a slightly different form, was suggested to me by Sam Gitler. In Section \ref{topsec}, we present our understanding of the
algebraic topology which motivated it.

\section{Enumeration when $k=4$}
In this section, we derive the 50 equivalence classes of configurations of size 4. We will abbreviate $\{i,j,k\}$ as $ijk$, etc. We will make constant use of the
fact that if $ij\sim  ijk$, then $ij\ell\sim ijk\ell$, which we sometimes call {\it unioning}. We will not  record equivalences of sets of the same size, since they are always a consequence
of other equivalences. For example, if $12\sim13$, then both are equivalent to $123$, and if $12\sim34$, then both are equivalent to $1234$.
We divide into cases according to the number of 3-sets which are equivalent to $1234$. We think of relations $\sim1234$ as being ``at the top,''
while relations $ij\sim ijk$ which are not $\sim1234$ as being ``at the bottom.''

{\bf Case 0}: No 3-sets equivalent to $1234$. Then there can be no equivalences at all (except that a set is equivalent to itself), due to unioning.
Thus there is one configuration of this type.

{\bf Case 1}: Exactly one 3-set $\sim1234$. WLOG $123\sim1234$. We cannot have, say, $12\sim123\sim1234$, since this would imply $124\sim1234$, which does not hold.
 It can be the case that either 0, 1, 2, or 3 of the relations $12\sim124$, $13\sim134$, and $23\sim234$ hold. (Note that these can be thought of as being obtained by dividing the relation $123\sim1234$ by 3 or 2 or 1.)
Any other $ij\sim ijk$ would imply a 3-set $\sim1234$ which is not the case. Thus there are four equivalence classes of configurations of this type.

{\bf Case 2}: Exactly two 3-sets $\sim1234$. WLOG these are 123 and 124. The only 2-set that might $\sim1234$ is 12. Thus there are two possibilities for what $\sim1234$,
depending on whether or not 12 is included.  This choice does not affect the possibilities for equivalences $ij\sim ijk$ which are not $\sim1234$. These possibilities are
$23\sim234$, $24\sim234$, $13\sim134$, and $14\sim134$, which may be thought of as dividing the top relations by appropriate integers. Note that if exactly one of these
four relations holds, then it doesn't matter which one, since any two are related by possibly interchanging 1 and 2, or possibly interchanging 3 and 4, and these do not
affect the relations at the top. Similarly if three of the four relations hold, it doesn't matter which three; each situation is equivalent under a permutation.
So far we have found four possibilities at the bottom, obtained by choosing 0, 1, 3, or 4 of the four possible equivalences. There are also three possible equivalence classes in
which we choose two of the four possible equivalences at the bottom. These are obtained by choosing $23\sim234$ and any one of the other three.
It is easily verified that no two of these three possibilities are equivalent under interchanges of 1 and 2 or of 3 and 4, which are the permutations that leave the top part fixed. Thus there are $2(4+3)=14$ equivalence classes of configurations
of this type.

{\bf Case 3}: Exactly three 3-sets $\sim1234$. These may be assumed to be 123, 124, and 134. We may also have 0, 1, 2, or 3 of 12, 13, and 14 equivalent to 1234. If none of them or all of them are $\sim1234$, then we obtain unique equivalence classes of configurations in which 0, 1, 2, or 3 of 23, 24, and 34 are equivalent to 234. If, say, we have just $12\sim 1234$ (and not 13 or 14), then 2 has a different status than 3 and 4 (which are interchangeable).  In this case, we can have 0 or 3 of 23, 24, and 34 equivalent to 234, or just 23, or just 34, or 23 and 24, or 23 and 34, so six possibilities at the bottom.
A similar situation occurs if two of the three (12, 13, and 14) are equivalent to 1234; one of 2, 3, and 4 will have a different status than the others, and so there will be six possibilities at the bottom. Thus there are $2\cdot4+2\cdot6=20$ equivalence classes of configurations of this type.

{\bf Case 4}: All four 3-sets $\sim1234$. Note that in this case it is impossible to have a relation $ij\sim ijk$ at the bottom, since $ijk\sim1234$. Of the $\binom42=6$ 2-sets, we can have any number, $t$, from 0 to 6 of them $\sim1234$. For $t=0$, 1, 5, or 6, there is just one
equivalence class of such configuration. When $t=2$, the configuration when $12\sim1234$ and $23\sim1234$ is not equivalent to the one when $12\sim1234$ and $34\sim1234$.
Similarly when $t=4$, there are two inequivalent configurations. When $t=3$, there are three inequivalent configurations,
represented by having the three sets $\sim 1234$ be $\{12,13,14\}$, $\{12,23,34\}$, or $\{12,13,23\}$.
 Thus there are $4\cdot1+2\cdot2+3=11$ configurations of this type.

Thus there are $1+4+14+20+11=50$ equivalence classes of configurations altogether.

\section{Tabulation and realization}
In Tables \ref{t2} and \ref{t3}, we list representatives of the 50 equivalence classes of configurations of size 4 obtained above.
Table \ref{t2} handles Cases 0 to 3 above, while Table \ref{t3} handles Case 4.
We use a dash instead of $\sim$ to improve the spacing in the table. We also list explicit small sets realizing each configuration.
For example, looking at the third entry in Table \ref{t2}, if $X_1=\{a,d\}$, $X_2=\{b,e\}$, $X_3=\{c\}$, and $X_4=\{d,e\}$, then $X_1\cup X_2=X_1\cup X_2\cup X_4$,
and no other equalities hold, except for that obtained by unioning with $X_3$.

\begin{table}[H]
\caption{$k=4$, Cases 0 to 3}\label{t2}
\begin{tabular}{cc|cccc}
Sets equaling 1234&Other equivalences&$X_1$&$X_2$&$X_3$&$X_4$\\
\hline
$\emptyset$&$\emptyset$&a&b&c&d\\
123&$\emptyset$&ad&be&cf&ef\\
123&12-124&ad&be&c&de\\
123&12-124,\,13-134&ad&be&ce&de\\
123&12-124,\,13-134,\,23-234&ade&bdf&cef&def\\
123,\,124&$\emptyset$&ace&bdf&cdg&efg\\
123,\,124&23-234&ac&bde&cdf&ef\\
123,\,124&23-234,\,24-234&a&bcd&ce&de\\
123,\,124&23-234,\,13-134&ace&bde&cdf&ef\\
123,\,124&23-234,\,14-134&ac&bd&ce&de\\
123,\,124&23-234,\,24-234,\,13-134&ad&bcd&ce&de\\
123,\,124&23-234,\,24-234,\,13-134,\,14-134&acd&bcd&ce&de\\
12,\,123,\,124&$\emptyset$&ace&bdf&cd&ef\\
12,\,123,\,124&23-234&acf&bde&cdf&ef\\
12,\,123,\,124&23-234,\,24-234&ae&bcd&ce&de\\
12,\,123,\,124&23-234,\,13-134&acef&bdeg&cdfg&efg\\
12,\,123,\,124&23-234,\,14-134&ace&bdf&cef&def\\
12,\,123,\,124&23-234,\,24-234,\,13-134&ade&bcdf&cef&def\\
12,\,123,\,124&23-234,\,24-234,\,13-134,\,14-134&acde&bcdf&cef&def\\
123,\,124,\,134&$\emptyset$&abce&bdf&cdg&efg\\
123,\,124,\,134&23-234&abc&bde&cdf&ef\\
123,\,124,\,134&23-234,\,24-234&ab&bcd&ce&de\\
123,\,124,\,134&23-234,\,24-234,\,34-234&a&bc&bd&cd\\
12,\,123,\,124,\,134&$\emptyset$&abce&bdf&cd&ef\\
12,\,123,\,124,\,134&23-234&abcf&bde&cdf&ef\\
12,\,123,\,124,\,134&34-234&abd&ce&bc&de\\
12,\,123,\,124,\,134&23-234,\,24-234&abe&bcd&ce&de\\
12,\,123,\,124,\,134&23-234,\,34-234&abe&cd&bce&de\\
12,\,123,\,124,\,134&23-234,\,24-234,\,34-234&ad&bc&bd&cd\\
13,\,14,\,123,\,124,\,134&$\emptyset$&abcd&bf&cef&def\\
13,\,14,\,123,\,124,\,134&23-234&abcd&bdf&cef&def\\
13,\,14,\,123,\,124,\,134&34-234&abcde&beg&bcfg&defg\\
13,\,14,\,123,\,124,\,134&23-234,\,24-234&abcd&bcdf&cef&def\\
13,\,14,\,123,\,124,\,134&23-234,\,34-234&abcd&bdf&bcef&def\\
13,\,14,\,123,\,124,\,134&23-234,\,24-234,\,34-234&abc&bce&bde&cde\\
12,\,13,\,14,\,123,\,124,\,134&$\emptyset$&abcd&be&ce&de\\
12,\,13,\,14,\,123,\,124,\,134&23-234&abcde&bdf&cef&def\\
12,\,13,\,14,\,123,\,124,\,134&23-234,\,24-234&abcde&bcdf&cef&def\\
12,\,13,\,14,\,123,\,124,\,134&23-234,\,24-234,\,34-234&abcd&bce&bde&cde
\end{tabular}
\end{table}

\begin{table}[H]
\caption{$k=4$, All 3-sets $\sim1234$}\label{t3}
\begin{tabular}{c|cccc}
2-sets equaling 1234&$X_1$&$X_2$&$X_3$&$X_4$\\
\hline
$\emptyset$&abd&ace&bcf&def\\
12&abd&ace&bc&de\\
12,\,23&abe&acde&bc&de\\
12,\,34&ac&bd&ab&cd\\
12,\,13,\,14&abc&ad&bd&cd\\
12,\,23,\,34&ade&bcd&abe&cde\\
12,\,13,\,23&abde&acdf&bcef&def\\
12,\,13,\,24,\,34&bcde&abdf&abef&cdef\\
12,\,13,\,14,\,24&abcd&abe&bde&cde\\
12,\,13,\,14,\,23,\,24&abcd&abce&bde&cde\\
12,\,13,\,14,\,23,\,24,\,34&abc&abd&acd&bcd
\end{tabular}
\end{table}

Sets realizing a certain configuration can be obtained systematically, although we obtained those in the tables using
a {\tt Maple} program.  We use the 17th entry in Table \ref{t2} to illustrate how the sets can be obtained
without using a computer. We want sets $X_i$ such that all relations are implied by $X_1\cup X_2=X_1\cup X_2\cup X_3\cup X_4$,
$X_2\cup X_3=X_2\cup X_3\cup X_4$, and $X_1\cup X_4=X_1\cup X_3\cup X_4$. In terms of inclusions, this says $X_3\subset X_1\cup X_2$, $X_4\subset X_1\cup X_2$, $X_4\subset X_2\cup X_3$,
$X_3\subset X_1\cup X_4$, $X_4\not\subset X_1\cup X_3$, $X_3\not\subset X_2\cup X_4$, $X_2\not\subset X_1\cup X_3\cup X_4$, and $X_1\not\subset X_2\cup X_3\cup X_4$. The ``proper'' condition adds the requirements $X_3\not\subset X_1$ and $X_4\not\subset X_2$.

Table \ref{t4} serves as a Venn diagram. Columns labeled 1 or 2 refer to just $X_1$ or $X_2$ without the other, and the first column means ``Neither $X_1$ nor $X_2$,'' with similar notation
for the rows. The condition $X_3\subset X_1\cup X_2$ forces the $\emptyset$ entries in the second and third entries in the first column, and the other $\emptyset$ entries
are forced similarly by the other inclusions above. The six noninclusions force, respectively, the entries $d$, $c$, $b$, $a$, $f$, and $e$ in Table \ref{t4}.

\begin{table}[H]
\caption{Venn diagram for one configuration}\label{t4}
\begin{tabular}{c|c|c|c|c|}
\mc{2}{c}{}&\mc{1}{c}{$1$}&\mc{1}{c}{$1\&2$}&\mc{1}{c}{$2$}\\
\cline{2-5}
&&$a$&&$b$\\
\cline{2-5}
$3$&$\emptyset$&$c$&&$\emptyset$\\
\cline{2-5}
$3\&4$&$\emptyset$&$e$&&$f$\\
\cline{2-5}
$4$&$\emptyset$&$\emptyset$&&$d$\\
\cline{2-5}
\end{tabular}
\end{table}

For future work, we might wish to know which pairs of sets have nontrivial intersection. In this example, from Table \ref{t4}, all pairs intersect except $X_1$ and  $X_2$.
This much would be forced. Having $X_1\cap X_2\ne\emptyset$ could be obtained by placing an element in any of the blank spaces in Table \ref{t4}.

\section{Topological motivation}\label{topsec}
First of all, our enumeration problem for sets is an interpretation of an enumeration problem for square-free monomials.
Let $m_1,\ldots,m_k$ be square-free monomials in a set of variables, such that no $m_i$ divides $m_j$ for $i\ne j$.
If $S\subset[k]$, let $m_S$ denote the least common multiple (lcm) of those $m_i$ for which $i\in S$. Define a lattice on $2^{[k]}$
by $S\le T$ if $m_S$ divides $m_T$. Since lcm is  the product  of the variables involved in the union of a set of monomials, this is clearly equivalent to our lattice $L(\M)$
defined at the outset.

Let $R$ be a commutative ring with 1.
Given a simplicial complex $K$ with vertex set $[n]$, there is an ideal $I(K)$ in $R[x_1,\ldots,x_n]$ generated by all monomials $x_{i_1}\cdots x_{i_r}$ such that $(i_1,\ldots,i_r)\not\in K$.
(For future reference, the quotient $R[x_1,\ldots,x_n]/I(K)$ is called the Stanley-Reisner ring $R(K)$.)
If we choose a minimal set of monomials generating this ideal, then this set of monomials will satisfy the condition that no $m_i$ divides $m_j$. Conversely, given a set of square-free
monomials $M$ in variables $x_1,\ldots,x_n$, there is a simplicial complex $K$ with vertex set $[n]$ such that $\sigma\not\in K$
iff some subset $\{i_1,\ldots,i_r\}$ of the vertices of $\sigma$ satisfies that $x_{i_1}\cdots x_{i_r}\in M$.

A simplicial complex $K$ with vertex set $[n]$ gives rise to a topological space called a moment-angle complex $Z_K$ defined in \cite[p.88]{BP}. In \cite[p.103]{BP}, it is proved that
there is an isomorphism
\begin{equation}\label{SR}H^*(Z_K;R)\approx \text{Tor}^{*,*}_{R[x_1,\ldots,x_n]}(R(K),R),\end{equation}
where grading in the RHS is by total degree.
Using the  Taylor resolution (e.g., \cite[p.439]{Eis}), these Tor groups are  related to our lattices $L(\M)$.

If the ideal $I(K)$ is spanned by monomials $m_1,\ldots,m_k$ in variables $x_1,\ldots,x_n$, the Taylor resolution leads to a cochain complex $C(K)$ whose cohomology is isomorphic to the RHS of (\ref{SR}). In grading $-j$, $C(K)$ is a free $R$-module on classes $e_S$ for all $S\subset [k]$ of cardinality $j$.
 The boundary homomorphism sends  $e_S$ to an alternating sum
of classes $e_{S-\{j\}}$ for which $\M_{\{j\}}\subset \M_{S-\{j\}}$. From the point of view of our tables, $e_S$ maps to an alternating sum of classes $e_T$ for which $T$ has cardinality $j-1$
and $T\sim S$.

The second grading in the Tor group is determined by giving $e_S$ grading equal to twice the cardinality of $\M_S$ in (\ref{Meq}). From the point of view of monomials, $|e_S|=\deg(\lcm\{m_i:i\in S\})$. Thus the isomorphism (\ref{SR}) depends not just on the lattice $L(\M)$ but also on the sets $X_i$ which led to the lattice. For example, the first lattice in Table \ref{t1} could have been realized either by $(X_1,X_2,X_3)=(\{a\},\{b\},\{c\})$ or $(\{a,b,c\},\{c,d,e\},\{e,f,a\})$. These give rise to different simplicial complexes,
different Stanley-Reisner rings, different moment-angle complexes, and different bigraded Tor groups, even though the lattices are equal. Their Tor groups are isomorphic if the second
bigrading is ignored, but this bigrading is essential to the isomorphism (\ref{SR}).

The isomorphism (\ref{SR}) is in fact an isomorphism of rings. The ring structure of the RHS depends on which pairs of monomials have gcd $>1$, or equivalently by which sets intersect.
The question posed to me by Sam Gitler asked to take into account in the enumeration also this information about overlaps. This can be done, but is more complicated.
We have chosen to focus here on the simpler question just based on inclusions.
\def\line{\rule{.6in}{.6pt}}

\end{document}